\documentclass{amsart}
\usepackage{palatino}
\usepackage{amsmath}
\usepackage{mathrsfs}
\usepackage[hmargin=3cm, vmargin=3cm]{geometry}
\usepackage[breaklinks,bookmarksopen,bookmarksnumbered]{hyperref}
\usepackage[all]{xy}

\theoremstyle{plain}
\newtheorem*{prop}{Proposition}

\newcommand\Q{\mathbb{Q}}
\renewcommand\P{\mathbb{P}}
\renewcommand\O{\mathscr{O}}

\begin{document}
\title{A non-hyperelliptic curve with torsion Ceresa class}
\author[Arnaud Beauville]{Arnaud Beauville}
\address{Universit\'e C\^ote d'Azur\\
CNRS -- Laboratoire J.-A. Dieudonn\'e\\
Parc Valrose\\
F-06108 Nice cedex 2, France}
\email{arnaud.beauville@unice.fr}
 
\begin{abstract}
We exhibit a non-hyperelliptic curve $C$ of genus 3 such that the class of the Ceresa cycle $[C]-[-C]$ in the intermediate Jacobian of $JC$ is torsion. 
\end{abstract}
\maketitle 
\section{Introduction}
Let $C$ be a complex curve  of genus $g\geq 3$, and $p$ a point of $C$. We embed $C$ into its Jacobian $J$ by the Abel-Jacobi map $x\mapsto [x]-[p]$. The 
 \emph{Ceresa cycle} $\mathfrak{z}_p(C)$ is the cycle $[C] - [(-1_J)^*C]$ in the Chow group $CH_1(J)_{\operatorname{hom} }$ of homologically trivial 1-cycles. The \emph{Ceresa class} $\mathfrak{c}_p(C)$ is the image of $\mathfrak{z}_p(C)$ in the intermediate Jacobian $\mathfrak{J}_1(J)$ parameterizing 1-cycles under the
 Abel-Jacobi map  $CH_1(J)_{\operatorname{hom} }\rightarrow \mathfrak{J}_1(J)$.

When $C$ is general,  $\mathfrak{z}_p(C)$ is not algebraically trivial \cite{C}. On the other hand, if $C$ is hyperelliptic $\mathfrak{z}_p(C)$ is algebraically trivial -- in fact it is zero if one chooses for $p$ a Weierstrass point. Not much is known
 besides these two extreme cases. There are  few curves for which $\mathfrak{z}_p(C)$ is known to be not 
algebraically trivial: Fermat curves of degree $\leq 1000$ \cite{O}, and the Klein quartic \cite{T}. An essential ingredient of these results is the fact that $\mathfrak{c}_p(C)$ is not a torsion class.

It is an open question whether there are non-hyperelliptic curves with $\mathfrak{z}_p(C)$ algebraically trivial. As observed in \cite[Remark 2.4]{FLV}, this condition is equivalent to a number of interesting 
properties: in particular the existence of a \emph{multiplicative Chow-K\"unneth decomposition} modulo algebraic equivalence, or the fact that the class $[C]\in CH_1(J)\otimes \Q$ is algebraically equivalent  to the minimal class $\dfrac{\theta ^{g-1}}{(g-1)!} $, where $\theta \in CH^1(J)$ is the class of the principal polarization.

In this note we exhibit a curve $C$ of genus 3 with the weaker property that the Ceresa class $\mathfrak{c}_p(C)$ is torsion (under the Bloch-Beilinson conjectures, this actually implies the algebraic triviality of $\mathfrak{z}_p(C)$ up to torsion). The construction is very simple: the curve $C$ has an automorphism $\sigma $ which fixes a point $p$, and therefore preserves $\mathfrak{c}_p(C)$; we just have to check that the fixed point set of $\sigma $ acting on $\mathfrak{J}_1(J)$ is finite.

 An analogous example, based on a much more sophisticated approach, appears in \cite{Litt}. It is not clear to me whether the Ceresa class in our sense of this example   is torsion.
 
 \bigskip	
 \section{The result}
\begin{prop}
Let $C\subset \P^2$ be the genus $3$ curve defined by $X^4+XZ^3+Y^3Z=0$, and let $p=(0,0,1)$. The Ceresa class  $\mathfrak{c}_p(C)$ is torsion. 
\end{prop}
\noindent\textit{Proof} :  Let $\omega  $ be a primitive 9-th root of unity. We consider the automorphism $\sigma $ of $C$ defined by $\sigma (X,Y,Z)=(X,\omega ^2Y,\omega ^3Z)$. We have $\sigma (p)=p$; therefore $\sigma $ preserves the Ceresa cycle $\mathfrak{z}_p(C)$, and also its class $\mathfrak{c}_p(C)$ in $\mathfrak{J}:=\mathfrak{J}_1(J)$. 

Thus it suffices to prove that $\sigma $ has finitely many fixed points on $\mathfrak{J}$;  equivalently, that the eigenvalues of $\sigma $ acting on  the tangent space $T_0(\mathfrak{J})$ are $\neq 1$. 

Now $T_0(\mathfrak{J})$ is identified with  $H^{0,3}(J)\,\oplus\, H^{1,2}(J)=\bigwedge^3V^*\oplus (\bigwedge^2V^*\otimes V)$, where
$V=H^{1,0}(J)=\allowbreak H^0(C,K_C)$. We first compute the eigenvalues of $\sigma $ on $V$. The elements of $V$ are of the form  $L\cdot \dfrac{XdZ-ZdX}{Y^2Z} $, with $L\in H^0(\P^2,\O_{\P}(1))$; it follows that the eigenvalues of $\sigma $ on $V$ are $\omega ^5,\omega ^7,\omega ^8$. Therefore the eigenvalue on 
$\bigwedge^3V^*$ is  $\omega ^7$, and the eigenvalues on   $\bigwedge^2V^*$ are $\omega ^3,\omega ^5,\omega ^6$. Thus each product of an eigenvalue on 
$\bigwedge^2V^*$ and one on $V$ is $\neq 1$, hence the Proposition.\qed

\end{document}